# Tight Chromatic Number for $\{K_{1,3}, (K_2 \cup K_1)+K_2\}$-free Graphs

Medha Dhurandhar

**ABSTRACT**

Although the chromatic number of a graph is not known in general, attempts have been made to find good bounds for the number. Here we prove that if a graph G is $\{K_{1,3}, (K_2 \cup K_1)+K_2\}$-free, and if $\Delta(G) \leq 2\omega(G)-3$, then $\chi(G)$ equals its maximum clique size. We also give examples of $\{K_{1,3}, (K_2 \cup K_1)+K_2\}$-free graphs with $\Delta(G) \geq 2\omega(G)-2$, where $\chi(G) = \omega(G)+1$.

## 1. INTRODUCTION

It has been an eminent unsolved problem in graph theory to determine the chromatic number of a given graph. Failing in the efforts to determine this, attempts have been made to find good bounds for the chromatic number of a graph. Vizing [4] proved that if a graph does not induce some nine subgraphs, then $\omega(G) \leq \chi(G) \leq \omega(G)+1$ where $\omega(G)$ is the size of maximum clique in G and $\chi(G)$ is the chromatic number of G. Later Choudum [1] and Javdekar [2] improved this result by dropping five and six of these nine subgraphs from the hypothesis, respectively. Finally Kierstead [3] showed that $\omega(G) \leq \chi(G) \leq \omega(G)+1$ for a $\{K_{1,3},(K_5-e)\}$-free graph. Moreover, Dhurandhar [5] proved that $\omega(G) \leq \chi(G) \leq \omega(G)+1$ for a $\{K_{1,3}, (K_2 \cup K_1)+K_2\}$-free graph. In this paper we prove that if G is $\{K_{1,3}, (K_2 \cup K_1)+K_2\}$-free, and if $\Delta(G) \leq 2\omega(G)-3$, then $\chi(G) = \omega(G)$.

**Notation:** For a graph G, V(G), E(G), $\Delta(G)$, $\omega(G)$, $\chi(G)$, and deg u denote the vertex set, edge set, maximum degree, size of a maximum clique, chromatic number, and the degree of u in G respectively. For a vertex $v \in V(G)$, $N(u) = \{v \in V(G) / uv \in E(G)\}$, and $\overline{N(u)} = N(u) \cup (u)$. If $S \subseteq V(G)$, then $<S>$ denotes the subgraph of G induced by S. If C is some colouring of G and if a vertex u of G is coloured m in C, then u is called a m-vertex. All graphs considered henceforth are simple.

Before proving the main result, we prove some lemmas.

**Lemma 1:** Let G be $\{K_{1,3}, (K_2 \cup K_1)+K_2\}$-free. Then $\forall\ u \in V(G)$, if $<Q>$ is a maximum clique in $<N(u)>$ and $R = N(u) - Q$, either
   1. $<N(u)> = C_5$ or
   2. $<N(u)> = P_4$ or
   3. $<R>$ is complete and every vertex of R is non-adjacent to a unique, distinct vertex of Q or
   4. $<R>$ is complete and $rq \notin E(G)\ \forall\ r \in R$ and $q \in Q$.

Proof:

**Case 1:** $|R| = 1$.

Clearly as G is $\{(K_2 \cup K_1)+K_2\}$-free, either r is non-adjacent to all vertices of Q or it is non-adjacent to only one vertex of Q.

**Case 2:** $|R| \geq 2$.

**Case 2.1:** $\exists\ x, y \in R$ s.t. $xy \notin E(G)$.
Let $z, w \in Q$ be s.t. $xz, yw \notin E(G)$. As G is $K_{1,3}$-free $z \neq w$ and $xw, yz \in E(G)$.

**Case 2.1.1:** $|N(u)| = 4$.
Then $<N(u)> = P_4$.

**Case 2.1.2:** $|N(u)| > 4$.

If $|Q| \geq 3$ and $q \in Q$, then as G is $\{(K_2 \cup K_1)+K_2\}$-free, clearly xq, yq $\in$ E(G). But then <x, w, q, u, y> = $(K_2 \cup K_1)+K_2$, a contradiction. Hence $|Q| = 2$ and $|R| \geq 3$. Let $v \in R$ and w.l.g. let vx $\in$ E(G). Now vw $\notin$ E(G) (else <x, v, w> is a bigger clique than Q in N(u)) $\Rightarrow$ vy $\in$ E(G). Also clearly $|R| = 2$ and <N(u)> = $C_5$.

**Case 2.2:** <R> is complete.

As $|Q| \geq |R|$, every vertex of R is non-adjacent to a distinct vertex of Q. Let $\exists$ x $\in$ R and y, z $\in$ Q s.t. xy, xz $\notin$ E(G). Clearly as G is $\{(K_2 \cup K_1)+K_2\}$-free, xq $\notin$ E(G) $\forall$ q $\in$ Q. Let if possible $\exists$ w $\in$ R s.t. say wy $\in$ E(G). W.l.g. let wz $\notin$ E(G). Thus either <N(u)> = $P_4$ or if <N(u)> $\neq P_4$, then clearly $|Q| \geq |R| \geq 3$. Let q $\in$ Q. Now if qw $\in$ E(G), then <y, q, u, w, x> = $(K_2 \cup K_1)+K_2$ and if qw $\notin$ E(G), then <q, z, u, y, w> = $(K_2 \cup K_1)+K_2$, a contradiction. Thus wv $\notin$ E(G) $\forall$ w $\in$ R and v $\in$ Q.

This proves Lemma 1.

**Corollary 1:** If G is $\{K_{1,3}, (K_2 \cup K_1)+K_2\}$-free, then $\Delta \leq 2\omega-1$ and either
  1. $\Delta = 2\omega-1$ with $\Delta = 5$, $\omega = 3$ or
  2. $\Delta \leq 2\omega-2$

For completeness we simply state some results from [5], which will be used in the main result.

**Result 1:** If G is $K_{1,3}$-free, then each component of the subgraph of G induced by two colour classes is either a path or a cycle.

**Result 2:** If G is $\{K_{1,3}, (K_2 \cup K_1)+K_2\}$-free, then $\omega(G) \leq \chi(G) \leq \omega(G) + 1$.

## 2. MAIN RESULT

**THEOREM**: If a graph G is $\{K_{1,3}, (K_2 \cup K_1)+K_2\}$-free, and if $\Delta(G) \leq 2\omega(G)-3$, then $\chi(G) = \omega(G)$.
Proof: Let if possible $\exists$ a $\{K_{1,3}, (K_2 \cup K_1)+K_2\}$-free graph G with $\Delta(G) \leq 2\omega(G)-3$ s.t. $\chi(G) > \omega(G)$. By **Result 2**, $\chi(G) = \omega(G)+1$. Let G be such a graph with minimum number of vertices.

Let u $\in$ V(G). If $\Delta(G-u) \leq 2\omega(G-u)-3$, then by minimality $\chi(G-u) = \omega(G-u) \geq \chi(G)-1 = \omega(G)$, and $\chi(G-u) = \omega(G)$. If $\Delta(G-u) > 2\omega(G-u)-3$, then $2\omega(G-u)-2 \leq \Delta(G-u) \leq \Delta(G) \leq 2\omega(G)-3 \Rightarrow \omega(G-u) = \omega(G)-1$ and by **Result 2**, $\chi(G-u) \leq \omega(G)$. Thus in any case $\chi(G-u) = \omega(G)$ $\forall$ u $\in$ V(G).

Let $C = \bigcup_{1}^{\omega+1} i$ be a $(\omega(G)+1)$-coloring of G s.t. u receives color $\omega(G)+1$.

**Case 1:** $\omega(G) \geq 4$.
Clearly $\forall$ v on a clique of size $\omega(G)$, <N(v)> $\neq C_5$, $P_4$. Hence by Lemma 1, we consider two subcases.

**Case 1.1:** $\exists$ u $\in$ V(G) lying on a clique Q of size $\omega$, with R = N(u) – Q s.t. every vertex of R is non-adjacent to a unique, distinct vertex of Q.
W.l.g. let u have maximum degree amongst all such vertices. Clearly every vertex of Q is non-adjacent to at the most one vertex of R. Let x $\in$ R have color $\omega$ in C s.t. $\omega$ is used in R but not in Q. Let y $\in$ Q be s.t. xy $\notin$ E(G). W.l.g. let y be a 1-vertex. Then R has no 1-vertex (else if z is a 1-vertex in R, then <Q-y+{x,z}> is a clique in G larger than <Q>). Now y (x) has a $\omega$-vertex (1-vertex) adjacent outside N(u) (else color y (x) by $\omega$ (1), and u by 1 ($\omega$)). Also as deg y, deg x $\leq$ deg u; y, x have no more vertices adjacent outside N(u). As $\Delta(G) \leq 2\omega(G)-3$, $\exists$ at least one more color say 2 used in Q but not in R. Let z be the 2-vertex in Q. Then z is the only 2-vertex adjacent to both x, y. Color. Again as deg z $\leq$ deg u; z has at the most one vertex adjacent outside N(u). Hence either y is the only 1-vertex of z or x is the only ω-vertex of z. Let α ∈ {1, ω} be that color. Color both x, y by 2, z by α and u by {1, ω}-α, a contradiction.

**Case 1.2:** ∀ u lying on a clique Q of size ω with R = N(u) – Q, rq ∉ E(G) ∀ r ∈ R and q ∈ Q.
We relabel vertices of Q as follows. Let u = $v_0$. As Δ(G) ≤ 2ω(G)-3, ∃ a colour, say 1 ∈ C, which is not used in N($v_0$)-Q. Now Q- $v_0$ must have a 1-vertex (else color $v_0$ by 1). Label this 1-vertex as $v_1$. Again ∃ a colour, say 2 ∈ C, which is not used in N($v_1$)-Q and Q- $v_1$ must have a 2-vertex (else color $v_1$ by 2, and $v_0$ by 1). Label this 2-vertex as $v_2$ and so on. We proceed similarly to get a maximal sequence S = $v_0$, $v_1$, ...., $v_k$ of vertices of Q s.t. $v_i$ is an i-vertex and i+1 is a color used in Q but not in N($v_i$)-Q. Clearly k ≥ 2. By maximality ∃ t < k s.t. t+1 is not used in N($v_k$)-Q. W.l.g. let d ∈ C be the color used in R but not in Q. Then each of $v_0$, $v_t$, $v_k$ has a unique d-vertex (else color that vertex by d and all the previous vertices $v_j$ in S by j+1). Also they all have distinct d-vertices (else we get **Case 1.1**). Consider the component T containing $v_{t+1}$ s.t. every vertex in T is colored either d or t+1. By Lemma 2, T is a path as $\deg_T v_{t+1}$ = 1. As G is $K_{1,3}$-free and $v_{t+1}$ is the only t+1 vertex of both $v_t$ and $v_k$ clearly either $v_t$ or $v_k$ is not adjacent to any vertex in T- $v_{t+1}$. If $v_t$ is not adjacent to any vertex in T-$v_{t+1}$, then interchange colors in T and color $v_j$ by j+1, (j ≤ t), a contradiction. If $v_k$ is not adjacent to any vertex in T- $v_{t+1}$, then consider the component K containing the d-vertex of $v_k$ and vertices having a color either d or t+1. As before clearly K is a path and $v_{t+1}$ is not adjacent to any vertex in K. Interchange colors in K, color $v_k$ by d and $v_j$ by j+1, (j < k), a contradiction.

**Case 2:** ω(G) ≤ 3.
If ω(G) = 2, and Δ(G) ≤ 1, then the result is trivially true. If ω(G) = 3, then Δ(G) ≤ 3 and 3 = ω(G) ≤ χ(G) ≤ Δ(G) ≤ 3 ⇒ χ(G) = ω(G), a contradiction.

This proves the theorem.

**Corollary 2:** If G is {$K_{1,3}$, ($K_2 \cup K_1$)+$K_2$}-free, then either
  3. Δ = 5 = 2ω-1, ω = 3 and χ(G) = ω+1. Here G has $W_6$ as an induced subgraph or
  4. Δ = 2ω-2 and χ(G) ≤ ω+1 or
  5. Δ ≤ 2ω-3 and χ(G) = ω.

**Examples to show that the condition Δ(G) ≤ 2ω(G)-3 is necessary.**

1. Let H = $C_{2n+1}$ (n>1) and V(H) = $\bigcup_{1}^{2n+1} u_i$. Construct G by replacing $u_{2i+1}$ by $K_m$ for 0 ≤ i ≤ n-1. Then G is {$K_{1,3}$, ($K_2 \cup K_1$)+$K_2$}-free, ω(G) = m+1, Δ(G) = 2m = 2ω(G)–2 and χ(G) = m+2 = ω(G)+1.

2. G = $W_6$. Here ω(G) = 3, Δ(G) = 5 = 2ω(G)–1 and χ(G) = 4 = ω(G)+1.